\documentclass[a4paper,12pt]{article}
\usepackage[margin=26mm]{geometry}
\usepackage{amsmath,amssymb,amsthm,mathtools}
\usepackage{xcolor}
\usepackage{tikz}
\usetikzlibrary{arrows.meta,decorations.markings,calc}
\usepackage[unicode,hidelinks]{hyperref}
\hypersetup{pdfauthor={},pdfsubject={Sharp digirth thresholds with bounded out-degree exceptions}}
\newtheorem{theorem}{Theorem}[section]
\newtheorem{lemma}[theorem]{Lemma}

\newtheorem{proposition}[theorem]{Proposition}
\newtheorem{claim}{Claim}
\newenvironment{claimproof}[1]{\begin{proof}[Proof of Claim~\ref{#1}]}{\end{proof}}
\DeclareMathOperator{\dist}{dist}

\newcommand{\len}[1]{\lVert #1\rVert}

\theoremstyle{definition}
\newtheorem{problem}{Problem}[section]

\newcommand{\proofcase}[2]{\par\medskip\noindent\textbf{Case #1.}\enspace #2\par\smallskip}

\tikzset{dot/.style={circle,fill=black,inner sep=1.6pt},
 patharc/.style={dashed,thick,postaction={decorate},decoration={markings,mark=at position .55 with {\arrow{Stealth}}}},
 realarc/.style={thick,-{Stealth[length=2mm]}},
 gate/.style={circle,draw,fill=white,inner sep=3pt}}

\title{Sharp digirth bounds for $C(k,k)$-subdivisions in digraphs}
\author{Hao Zhang\;\quad Panpan Cheng\\
	\small{School of Mathematics, Shandong University, Shandong, Jinan, 250000, China}
	}
\date{}
\hypersetup{pdftitle={Sharp digirth bounds for C(k,k)-subdivisions in digraphs }}
\begin{document}
\maketitle
\begin{abstract}
For an integer \(k\ge2\), let \(C(k,k)\) denote the oriented cycle formed by two internally vertex-disjoint directed paths of length \(k\) with the same initial and terminal vertices. We extend a problem of Picasarri-Arrieta and Rambaud [European J. Combin. 122 (2024)] on subdivisions of \(C(k,k)\) to digraphs with several vertices of out-degree one. Following the sharp result of Chen et al. [arXiv:2610.38574v1 (2026)] for at most one such vertex, we determine the exact digirth threshold when at most \(r\ge4\) such vertices are permitted. Specifically, for all integers \(k\ge2\) and \(r\ge4\), every digraph of minimum out-degree at least one, with at most \(r\) vertices of out-degree one and digirth at least \(r(k-1)+1\), contains a subdivision of \(C(k,k)\). 
\end{abstract}

\section{Introduction}

A subdivision of a digraph $F$ is obtained by replacing each arc of $F$ by a directed path with the same initial and terminal vertices, where the replacement paths have pairwise disjoint interiors and their internal vertices avoid $V(F)$. A central question is whether sufficiently large minimum out-degree guarantees a subdivision of a prescribed digraph.

Mader~\cite{mader1985} conjectured that, for every acyclic digraph $F$, there exists a positive integer $d$ such that every digraph of minimum out-degree at least $d$ contains a subdivision of $F$. Not every digraph has this property: Thomassen~\cite{thomassen1985} constructed digraphs of arbitrarily large minimum out-degree with no even directed cycle. Mader~\cite{mader1995} also proved that sufficiently large minimum out-degree forces two vertices joined by any prescribed number of internally vertex-disjoint directed paths. Requiring these paths to be long leads to the study of cycles with two blocks.

For integers $k,\ell\ge2$, let $C(k,\ell)$ be the union of two internally vertex-disjoint directed paths of lengths $k$ and $\ell$ with the same initial and terminal vertices. Aboulker et al.~\cite{aboulker} proved that every digraph of minimum out-degree at least $2k+2\ell-1$ contains a subdivision of $C(k,\ell)$. For $k\ge\ell\ge2$, Gishboliner, Steiner and Szab\'o~\cite{gss} reduced this sufficient out-degree bound to $k+3\ell-5$. Liu and Yoo~\cite{liuyoo} subsequently proved that minimum out-degree at least $k+\ell-1$ suffices, and that this bound is best possible.
Subdivisions of two-block cycles and related path systems have also been investigated under chromatic assumptions~\cite{cohencycles,bispindles,almniny}.

Digirth conditions can reduce the minimum out-degree required to force a subdivision. Let $g(D)$ denote the length of a shortest directed cycle in $D$. Picasarri-Arrieta and Rambaud~\cite{pr} proved that sufficiently large digirth and minimum out-degree at least two force a subdivision of $C(k,k)$, and asked for the least possible digirth threshold. Chen, Hou, Ma, Yin and Zhou~\cite{chen} obtained the upper bound $4k+2$; their sufficient condition also permits one vertex of out-degree one. More recently, Chen, Hou, Hou, Liu, Wang and Zhou~\cite{chenexact} determined the exact threshold to be $k+1$, including the case in which one vertex may have out-degree one. Their proof establishes a structural theorem on two directed cycles whose intersection is exactly one arc.

We consider the corresponding problem when several vertices may have out-degree one. Without a bound on their number, no digirth condition can force a subdivision of $C(k,k)$: directed cycles have arbitrarily large digirth and contain no such subdivision. This motivates a threshold depending on both $k$ and the number of permitted exceptions. Vertices of out-degree one also occur in other digirth estimates. For example, Shen~\cite{shen} bounded the diameter of a strongly connected digraph of order $n$, digirth $g$ and $t$ vertices of out-degree one by $n-g+t$, while Aharoni et al.~\cite{nonuniform} obtained a digirth bound in terms of the individual out-degrees.

For a digraph $D$, let $O(D)=\{v\in V(D):d_D^+(v)=1\}$.
We call the vertices of $O(D)$ exceptional. For integers $k\ge2$ and $r\ge0$, let $F_r(k)$ be the least positive integer $g$ such that every digraph $D$ satisfying $\delta^+(D)\ge1$, $|O(D)|\le r$ and $g(D)\ge g$ contains a subdivision of $C(k,k)$. Set $F_r(k)=+\infty$ if no such integer exists. The original question of Picasarri-Arrieta and Rambaud~\cite{pr} concerns $F_0(k)$. Recently, Chen et al.~\cite{chenexact} obtained
$F_0(k)=F_1(k)=k+1$. This leads to the question of determining $F_r(k)$ for $r\ge2$. Our main result answers this question for every $k\ge2$ and $r\ge4$.


\begin{theorem}\label{thm:many}
For all integers $k\ge2$ and $r\ge4$,
$F_r(k)=r(k-1)+1$.

\end{theorem}

In concrete terms, for $k\ge2$ and $r\ge4$, Theorem~\ref{thm:many} states that every digraph of minimum out-degree at least one, with at most $r$ vertices of out-degree one and digirth at least $r(k-1)+1$, contains a subdivision of $C(k,k)$. At digirth $r(k-1)$, this conclusion can fail even when exactly $r$ vertices have out-degree one and every other vertex has out-degree two. Thus the digirth threshold is sharp.

When no exceptions are permitted, the degree condition reduces to minimum out-degree at least two, as in~\cite{pr,chen}.
For fixed $k$, the threshold therefore increases by exactly $k-1$ with each additional permitted exception once $r\ge4$.

The proof builds on the isometric-cycle method of Chen et al.~\cite{chen} and a local form of the minimal-set argument used in~\cite{chenexact}. The main additional step is to control several exceptional vertices simultaneously. In a minimal counterexample, we associate regions with short segments of a shortest directed cycle. Each region contains an exceptional vertex, whereas the union of any two intersecting regions contains every exceptional vertex. Together with restrictions on which regions can intersect, these properties force the regions to be pairwise disjoint. Counting them then yields the upper bound. A cyclic construction establishes sharpness.

Section~\ref{sec:prelim} introduces the notation and proves the local structural lemma. Section~\ref{sec:many} proves Theorem~\ref{thm:many}, including the extremal construction. Section~\ref{sec:remarks} discusses the remaining cases.

\section{Preliminaries}\label{sec:prelim}
All digraphs considered in this paper are finite and nonempty, have no loops or parallel arcs, and may contain pairs of opposite arcs. For a digraph $D=(V,A)$, its order and size are denoted by $|D|=|V(D)|$ and $\len D=|A(D)|$, respectively. We write $v\in D$ for $v\in V(D)$ and $uv\in D$ for $uv\in A(D)$. The out-neighborhood of $v$ is $N_D^+(v)=\{u\in V(D):vu\in A(D)\}$, its out-degree is $d_D^+(v)=|N_D^+(v)|$, and $\delta^+(D)=\min_{v\in V(D)}d_D^+(v)$. For $S\subseteq V(D)$, let $D[S]$ be the induced subdigraph. We write $D-v$ for vertex deletion and $D-uv$ for deletion of the arc $uv$ only.

All paths and cycles are directed and simple, and their lengths are their numbers of arcs. Single-vertex paths of length zero are allowed. The digirth $g(D)$ is the length of a shortest directed cycle, with $g(D)=+\infty$ if $D$ is acyclic. Every finite digraph of minimum out-degree at least one contains a directed cycle, so its digirth is finite. A path with initial vertex $u$ and terminal vertex $v$ is called a $(u,v)$-path. Its minimum possible length is denoted by $\dist_D(u,v)$, with value $+\infty$ if no such path exists.

For vertices $a,b$ of a cycle $C$, let $C[a,b]$ be the path from $a$ to $b$ along the direction of $C$, interpreted as a single-vertex path if $a=b$. The notation $C[a,b)$, $C(a,b]$ and $C(a,b)$ indicates deletion of the corresponding endpoints. The same notation is used for subpaths of a path. Membership in a path or cycle refers to membership in its vertex set. When two paths meet only at the terminal vertex of the first and the initial vertex of the second, their concatenation is denoted by $P\circ Q$. Concatenations that may repeat vertices will explicitly be called walks. Every directed walk contains a path with the same endpoints using only vertices of that walk.

For distinct vertices $u,v$ with $v$ reachable from $u$, a vertex $c\notin\{u,v\}$ is a \emph{$(u,v)$-cut-vertex} if every $(u,v)$-path contains $c$. We use the directed vertex form of Menger's theorem: if $u\ne v$, $uv\notin A(D)$, and $v$ is reachable from $u$, then either there are two internally vertex-disjoint $(u,v)$-paths or there is a $(u,v)$-cut-vertex. A cycle $C$ is \emph{isometric} in $D$ if $\dist_D(a,b)=\len{C[a,b]}$ for all $a,b\in C$. Every shortest directed cycle is isometric. Indeed, a shorter replacement for a segment of such a cycle, followed by the complementary segment, would give a nonempty closed directed walk shorter than the cycle. This walk contains a directed cycle, contradicting the choice of the original cycle.

For $S\subseteq V(D)$, define its external out-neighborhood by
\[
\partial_D^+(S)=\{w\in V(D)\setminus S:vw\in A(D)\text{ for some }v\in S\}.
\]
If $D$ is strongly connected, then $\partial_D^+(S)\ne\varnothing$ for every nonempty proper subset $S$ of $V(D)$. In particular, if $S$ consists of all vertices reachable from a vertex $p\ne q$ in $D-q$, then $\partial_D^+(S)=\{q\}$: any other out-neighbor of a vertex of $S$ would also be reachable from $p$ in $D-q$.

The following lemma localizes the minimal-set argument in the proof of \cite[Theorem~1.5]{chenexact}. Only the vertices of the prescribed set $X$ are required to have out-degree at least two, so the lemma applies in the presence of several exceptional vertices outside $X$. We include a proof for completeness. Its structural conclusion does not require a digirth assumption.

\begin{lemma}\label{lem:oneexit}
Let $D$ be a strongly connected digraph with $\delta^+(D)\ge1$. Suppose that $X$ is a nonempty proper subset of $V(D)$ satisfying $|\partial_D^+(X)|=1$ and $X\cap O(D)=\varnothing$. Then $D$ contains an arc $xy$ and two internally vertex-disjoint $(y,x)$-paths. Consequently, if $g(D)\ge k+1$ and $D$ contains no subdivision of $C(k,k)$, then every nonempty proper subset $X$ with $|\partial_D^+(X)|=1$ satisfies $X\cap O(D)\ne\varnothing$.
\end{lemma}

\begin{proof}
Choose a nonempty subset $S\subseteq X$ of minimum cardinality subject to $|\partial_D^+(S)|=1$, and write $\partial_D^+(S)=\{q\}$. Such a choice is possible since $X$ itself has this property. By strong connectivity, there is a path $R$ from $q$ to $S$ that first meets $S$ at its terminal vertex $x$. Hence $R$ has no internal vertex in $S\cup\{q\}$. Since $S\cap O(D)=\varnothing$, we have $d_D^+(x)\ge2$. The only possible out-neighbor of $x$ outside $S$ is $q$, so there is an arc $xy$ with $y\in S$. Put $H=D[S\cup\{q\}]$.

We claim that $H$ contains a $(y,x)$-path $P_x$ and a $(y,q)$-path $P_q$ such that $V(P_x)\cap V(P_q)=\{y\}$. To see this, add a new vertex $t$ to $H$ and add precisely the arcs $xt$ and $qt$. The desired paths exist if and only if the resulting digraph contains two internally vertex-disjoint $(y,t)$-paths. Indeed, two such paths must enter $t$ through distinct vertices of $\{x,q\}$, and deleting their final arcs produces the required pair.

Suppose that no such pair of paths exists. Since $y\notin\{x,q\}$, the arc $yt$ is absent. If $t$ is reachable from $y$, Menger's theorem gives a vertex $z\notin\{y,t\}$ whose deletion destroys every path from $y$ to $t$; set $Z=\{z\}$. Otherwise, set $Z=\varnothing$. In either case, $Z\subseteq V(H)\setminus\{y\}$, $|Z|\le1$, and there is no path from $y$ to $t$ after deleting $Z$. Let
\[
T=\{v\in V(H)\setminus Z:\text{there is a $(y,v)$-path in }H-Z\}.
\]
Neither $x$ nor $q$ belongs to $T$, since reaching either undeleted vertex would also allow a path to $t$. Thus $\varnothing\ne T\subseteq S\setminus\{x\}$. All out-neighbors in $D$ of vertices of $S$ lie in $S\cup\{q\}$. If $u\in T$ and $uv\in A(D)$ with $v\notin Z$, then $v$ is also reachable from $y$ in $H-Z$. Therefore $\partial_D^+(T)\subseteq Z$. Strong connectivity implies that $\partial_D^+(T)$ is nonempty, and hence $|\partial_D^+(T)|=1$. This contradicts the minimality of $S$, because $T\subseteq X$ and $|T|<|S|$.

Thus the required paths $P_x$ and $P_q$ exist. Since they lie in $H$ and meet only at $y$, while the internal vertices of $R$ lie outside $H$, the paths $P_x$ and $P_q\circ R$ are internally vertex-disjoint paths from $y$ to $x$. This proves the first assertion. Adding the arc $xy$ to either path forms a directed cycle. Hence, if $g(D)\ge k+1$, both paths have length at least $g(D)-1\ge k$, so their union is a subdivision of $C(k,k)$.
\end{proof}

\section{Proof of Theorem~\ref{thm:many}}\label{sec:many}
\subsection{\texorpdfstring{The upper bound for $k\ge3$}{The upper bound for k >= 3}}\label{sec:upper}
We first prove the upper bound in Theorem~\ref{thm:many}. Fix $k\ge3$ and $r\ge4$, and put $G=r(k-1)+1$. We have
\begin{equation}\label{eq:digirth}
G\ge4k-3\ge3k.
\end{equation}
The case $k=2$ is treated in Section~\ref{sec:ktwo}, and the construction establishing the lower bound is given in Section~\ref{sec:lower}.

Suppose that the assertion is false, and choose a counterexample $D$ with $|D|$ minimum and, subject to this, $\len D$ minimum. A terminal strongly connected component retains all outgoing arcs of its vertices and inherits both the digirth bound and the absence of the required subdivision. The minimality of $|D|$ therefore implies that $D$ is strongly connected. If a vertex had out-degree at least three, deleting one of its outgoing arcs would preserve all the conditions defining a counterexample and contradict the minimality of $\len D$. Thus every vertex has out-degree one or two. Moreover, $|D|\ge g(D)\ge G>r$. In particular, not every vertex can have out-degree one: otherwise $|D|=|O(D)|\le r$, and hence $g(D)\le|D|\le r$, contrary to $g(D)\ge G>r$. Thus a vertex of out-degree two exists. If $|O(D)|<r$, deleting one outgoing arc of a vertex of out-degree two would produce a counterexample with the same order, fewer arcs and at most $r$ vertices of out-degree one. Write $O=O(D)$. It follows that
\begin{equation}\label{eq:minimal}
d_D^+(v)\in\{1,2\}\;(v\in V(D))\;\text{and}\; |O|=r.
\end{equation}

Choose a shortest directed cycle $C$ in $D$ and write $L=\len C=g(D)$. Thus $C$ is isometric and $L\ge3k$. We will repeatedly use the following observation. If $x,y$ are distinct vertices of $C$ and $P$ is an $(x,y)$-path whose internal vertices avoid $C$, then
\begin{equation}\label{eq:ear}
\len{C[x,y]}\le k-1.
\end{equation}
Otherwise, isometry would give $\len P\ge\dist_D(x,y)=\len{C[x,y]}\ge k$, and $P\cup C[x,y]$ would be a subdivision of $C(k,k)$.

\begin{claim}\label{cl:chain}
For distinct $p,q\in C$, list all $(p,q)$-cut-vertices in their order on $C[p,q]$ as $c_1,\ldots,c_s$, and put $c_0=p$ and $c_{s+1}=q$. Every $(p,q)$-path visits these vertices in this order, and
\[
\dist_D(c_j,c_{j+1})=\len{C[c_j,c_{j+1}]}\le k-1
\;(0\le j\le s).
\]
\end{claim}

\begin{claimproof}{cl:chain}
Every $(p,q)$-cut-vertex lies internally on $C[p,q]$, since this segment is itself a $(p,q)$-path. Suppose that a $(p,q)$-path $P$ visits $c_j$ before $c_i$, where $i<j$. Then $P[p,c_j]$ followed by $C[c_j,q]$ is a $(p,q)$-walk avoiding $c_i$. It contains a $(p,q)$-path avoiding $c_i$, a contradiction. The asserted order follows.

There is no $(c_j,c_{j+1})$-cut-vertex. Indeed, such a vertex would lie internally on $C[c_j,c_{j+1}]$ and, by the order just established, would occur on every $(p,q)$-path. It would therefore be an additional $(p,q)$-cut-vertex between two consecutive members of the list. If $\dist_D(c_j,c_{j+1})\ge k$, there is no arc $c_jc_{j+1}$, and Menger's theorem gives two internally vertex-disjoint $(c_j,c_{j+1})$-paths. Both have length at least $k$, a contradiction. Isometry now implies the required bound. The argument also applies when the list of internal cut-vertices is empty.
\end{claimproof}

Let $uv$ be an arc of $C$. Adding $uv$ to any $(v,u)$-path produces a directed cycle, so every such path has length at least $L-1\ge k$. If there were no $(v,u)$-cut-vertex, Menger's theorem would give two internally vertex-disjoint $(v,u)$-paths, and hence a subdivision of $C(k,k)$, a contradiction. Now consider the first vertex in the cut-vertex sequence whose distance from $v$ is at least $k$. Its predecessor in the sequence has distance at most $k-1$ from $v$, and Claim~\ref{cl:chain} bounds the distance between consecutive vertices by $k-1$. Thus the chosen vertex has distance at most $2k-2$ from $v$. Since $\dist_D(v,u)=L-1>2k-2$, this vertex is an internal cut-vertex.

\begin{claim}\label{cl:normalize}
There is an arc $ab$ of $C$ for which the first $(b,a)$-cut-vertex at distance at least $k$ from $b$ has distance exactly $k$ from $b$.
\end{claim}

\begin{claimproof}{cl:normalize}
Start with an arbitrary arc $uv$ of $C$, and let $c$ be the first $(v,u)$-cut-vertex at distance at least $k$ from $v$. If $\dist_D(v,c)=k$, take $ab=uv$. Otherwise choose $b\in C[v,c]$ with $\len{C[b,c]}=k$, and let $a$ be the predecessor of $b$ on $C$.

We show that $c$ is a $(b,a)$-cut-vertex. Suppose that there is a $(b,a)$-path avoiding $c$. Before its first visit to $C[v,a]$, take its last vertex in $C[b,c)$, and denote the resulting subpath by $P$, with initial vertex $x\in C[b,c)$ and terminal vertex $y\in C[v,a]$. Its internal vertices avoid these two segments, and it avoids $c$. If an internal vertex $z$ of $P$ lies on $C$, then $z\in C(c,u]$. The walk $C[v,x]\,P[x,z]\,C[z,u]$ avoids $c$, contrary to its being a $(v,u)$-cut-vertex. Thus $P$ is internally disjoint from $C$. Since $\len{C[y,x]}\le\len{C[v,c]}\le2k-2$, we have $\len{C[x,y]}\ge L-(2k-2)\ge k+2$, contradicting~\eqref{eq:ear}. This proves the assertion about $c$. As all vertices preceding $c$ on $C[b,c]$ have distance less than $k$ from $b$, it is the required first cut-vertex.
\end{claimproof}

Fix an arc $ab$ as in Claim~\ref{cl:normalize}. Write all $(b,a)$-cut-vertices as $c_1,\ldots,c_\ell$ in their order on $C[b,a]$, with $c_0=b$ and $c_{\ell+1}=a$. Let $i_0$ be the least index with $\dist_D(b,c_{i_0})\ge k$, and let $i_1$ be the greatest index with $\dist_D(c_{i_1},a)\ge k$. Put $h=i_1-i_0$ and $\tilde c_i=c_{i_0+i}$ for $0\le i\le h$. Claims~\ref{cl:chain} and~\ref{cl:normalize} imply that
\begin{equation}\label{eq:coarse}
\begin{gathered}
\len{C[b,\tilde c_0]}=k,\;
\len{C[\tilde c_i,\tilde c_{i+1}]}\le k-1\; (0\le i<h)\;\text{and}\;
k\le\len{C[\tilde c_h,a]}\le2k-2.
\end{gathered}
\end{equation}
The indices are well defined: the remaining distance from $\tilde c_0$ to $a$ is $L-1-k\ge2k-1$. In particular $i_1\ge i_0$. If equality held, then~\eqref{eq:coarse} would give $L\le k+(2k-2)+1=3k-1$, a contradiction. Hence $h\ge1$.

In the following, indices on the vertices $\tilde c_i$ are read modulo $h+1$. Define
\[
X_i=\{v\in V(D):\text{there is a $(\tilde c_i,v)$-path in }D-\tilde c_{i+1}\}.
\]
We refer to $C[\tilde c_i,\tilde c_{i+1}]$, for $0\le i<h$, as the ordinary segments, and to the remaining segment $C[\tilde c_h,\tilde c_0]$, which passes through $a$ and $b$, as the special segment. These names distinguish the segments by their positions on $C$; they do not refer to the vertices in $O$.

\begin{claim}\label{cl:coarse}
For every $0\le i\le h$,
\[
X_i\cap V(C)=V(C[\tilde c_i,\tilde c_{i+1}))\;\text{and}\;\partial_D^+(X_i)=\{\tilde c_{i+1}\}.
\]
\end{claim}

\begin{claimproof}{cl:coarse}
The displayed half-open segment is contained in $X_i$. Suppose that a path in $D-\tilde c_{i+1}$ from $\tilde c_i$ reaches another vertex of $C$. Take its first vertex outside the designated segment on $C$ and, before this vertex, its last vertex in the designated segment. These vertices are the terminal and initial vertices, respectively, of a subpath $P$ internally disjoint from $C$. Write them as $y$ and $x$.

Suppose first that $i<h$. If $y\in C(\tilde c_{i+1},a]$, the walk $C[b,x]PC[y,a]$ avoids the $(b,a)$-cut-vertex $\tilde c_{i+1}$, a contradiction. Otherwise $y\in C[b,\tilde c_i)$. Then $C[x,y]$ contains $C[\tilde c_h,b]$, whose length is at least $k+1$ by~\eqref{eq:coarse}. This contradicts~\eqref{eq:ear}.

Now suppose that $i=h$. Then $y\in C(\tilde c_0,\tilde c_h)$. If $x\in C[b,\tilde c_0)$, the walk $C[b,x]PC[y,a]$ avoids $\tilde c_0$, again a contradiction. If $x\in C[\tilde c_h,a]$, then $C[x,y]$ contains $C[a,\tilde c_0]$, of length $k+1$, contradicting~\eqref{eq:ear}. This proves the equality for the intersection with $C$. The definition of $X_i$ implies that its external out-neighbors are contained in $\{\tilde c_{i+1}\}$, and strong connectivity makes this inclusion an equality.
\end{claimproof}

The special segment may have length greater than $k-1$. We therefore refine the partition by selecting additional vertices on each segment $C[p,q]$ as follows: retain its endpoints and select all $(p,q)$-cut-vertices in $D$. By Claim~\ref{cl:chain}, these vertices lie on $C[p,q]$, and each segment between consecutive selected vertices has length at most $k-1$. Counting each shared endpoint only once, we list the selected vertices in their cyclic order as $\tilde c_0^*,\tilde c_1^*,\ldots,\tilde c_{N-1}^*$.
Indices on starred vertices and sets will always be read modulo $N$. We have
\begin{equation}\label{eq:finegap}
1\le\len{C[\tilde c_j^*,\tilde c_{j+1}^*]}\le k-1
\;(0\le j<N).
\end{equation}
Summing~\eqref{eq:finegap} around $C$ yields $L\le N(k-1)$. Since $L\ge r(k-1)+1$, we obtain $N\ge r+1\ge5$. Define
\[
X_j^*=\{v\in V(D):\text{there is a $(\tilde c_j^*,v)$-path in }D-\tilde c_{j+1}^*\}.
\]

Figure~\ref{fig:partition} illustrates the refinement.

\begin{figure}[tbp]
\centering
\begin{tikzpicture}[font=\small]
\node[dot,label=below left:{$a$}] (a) at (-.85,-2.1) {};
\node[dot,label=below right:{$b$}] (b) at (.85,-2.1) {};
\node[dot,label=right:{$\tilde c_0$}] (t0) at (2.35,-.65) {};
\node[dot,label=above right:{$\tilde c_1$}] (t1) at (1.7,1.35) {};
\node[dot,label=above left:{$\tilde c_{h-1}$}] (tm) at (-1.7,1.35) {};
\node[dot,label=left:{$\tilde c_h$}] (th) at (-2.35,-.65) {};
\draw[realarc] (a) to[bend right=15] (b);
\draw[patharc] (b) to[bend right=24] node[right=5pt] {$k$} (t0);
\draw[patharc] (t0) to[bend right=20] node[right=5pt] {$\le k-1$} (t1);
\draw[patharc] (t1) to[bend right=30] node[above=5pt] {$\cdots$} (tm);
\draw[patharc] (tm) to[bend right=20] node[left=5pt] {$\le k-1$} (th);
\draw[patharc] (th) to[bend right=24] node[left=7pt] {$[k,2k-2]$} (a);
\node at (0,.55) {$C$};
\node at (0,-.1) {ordinary segments};
\node at (0,-1.0) {special segment};
\begin{scope}[yshift=-4.0cm]
\node[dot,label=below:{$p$}] (p) at (-3.5,0) {};
\node[dot,label=below:{$s=\tilde c_j^*$}] (ss) at (-1.2,0) {};
\node[dot,label=below:{$t=\tilde c_{j+1}^*$}] (tt) at (1.2,0) {};
\node[dot,label=below:{$q$}] (q) at (3.5,0) {};
\draw[patharc] (p)--(ss);
\draw[patharc] (ss)--node[above=4pt] {$\le k-1$}(tt);
\draw[patharc] (tt)--(q);
\node at (0,-.85) {Insert all $(p,q)$-cut-vertices in each segment $C[p,q]$};
\end{scope}
\end{tikzpicture}
\caption{The partition of $C$, illustrated for $h\ge3$. Solid lines denote arcs and dashed lines denote directed paths.}
\label{fig:partition}
\end{figure}

\begin{claim}\label{cl:fine}
For every $0\le j<N$,
\begin{equation}\label{eq:fineintersection}
X_j^*\cap V(C)=V(C[\tilde c_j^*,\tilde c_{j+1}^*))
\;\text{and}\; \partial_D^+(X_j^*)=\{\tilde c_{j+1}^*\}.
\end{equation}
Moreover, $X_j^*\cap O\ne\varnothing$.
\end{claim}

\begin{claimproof}{cl:fine}
Let $s=\tilde c_j^*$ and $t=\tilde c_{j+1}^*$, and let $C[p,q]$ be the segment before refinement whose cut-vertex sequence contains $s,t$ consecutively. Denote its associated set from Claim~\ref{cl:coarse} by $X$. We first show that $X_j^*\subseteq X$. If $t=q$, a path from $s$ avoiding $q$, preceded by $C[p,s]$, forms a walk from $p$ avoiding $q$. If $t\ne q$, then $t$ is a $(p,q)$-cut-vertex. There is no $(s,q)$-path avoiding $t$, since preceding such a path by $C[p,s]$ would contradict this property. As $s\in X$ and the only external out-neighbor of $X$ is $q$, no path from $s$ avoiding $t$ can leave $X$. The inclusion follows in both cases.

Consequently, all vertices of $X_j^*$ on $C$ belong to $C[p,q)$. None can lie after $t$ in the order on this segment. Indeed, if a path $P$ from $s$ to such a vertex $w$ avoids $t$, then $C[p,s]PC[w,q]$ is a $(p,q)$-walk avoiding $t$. This is impossible when $t\ne q$, and there are no such vertices when $t=q$.

It remains to show that no vertex of $C[p,s)$ belongs to $X_j^{*}$.
Suppose otherwise, and take a path in $D-t$ from $s$ to such a vertex.
Let $y$ be the first vertex of this path in $C[p,s)$, and let $x$ be
the last vertex before $y$ that lies in $C[s,t)$. Let $P$ be the
subpath from $x$ to $y$. Thus $x\in C[s,t)$ and $y\in C[p,s)$.
By the inclusion $X_j^{*}\subseteq X$ and the preceding paragraph,
every vertex of the chosen path that lies on $C$ belongs to $C[p,t)$.
Since $y$ is the first vertex in $C[p,s)$ and $x$ is the last vertex
in $C[s,t)$ before $y$, no internal vertex of $P$ belongs to either
segment. Hence all internal vertices of $P$ lie outside $C$.

\proofcase{1}{$C[p,q]$ is ordinary.}
Then $\len{C[y,x]}\le k-2$, since $x\ne q$. Thus $\len{C[x,y]}=L-\len{C[y,x]}\ge k$, contradicting~\eqref{eq:ear}.

\proofcase{2}{$C[p,q]=C[\tilde c_h,\tilde c_0]$.}
Its order runs through $\tilde c_h,\ldots,a,b,\ldots,\tilde c_0$. When $x,y\in C[\tilde c_h,a]$, we have $\len{C[y,x]}\le2k-2$ and therefore $\len{C[x,y]}\ge k+2$. When $x,y\in C[b,\tilde c_0)$, we have $\len{C[y,x]}\le k-1$ and again $\len{C[x,y]}\ge k$. Both possibilities contradict~\eqref{eq:ear}. The remaining possibility is $y\in C[\tilde c_h,a]$ and $x\in C[b,\tilde c_0)$. In this case $C[b,x]PC[y,a]$ avoids the $(b,a)$-cut-vertex $\tilde c_0$, a contradiction. These cases cover all possible positions of $x$ and $y$ on the special
segment when $y$ precedes $x$.

Every vertex of $C[s,t)$ can be reached from $s$ along this segment
without passing through $t$. Hence
$X_j^{*}\cap V(C)=V(C[s,t))$.
By the definition of $X_j^{*}$, every arc leaving $X_j^{*}$ ends at $t$.
Moreover, $X_j^{*}$ is nonempty and proper, since it contains $s$ but
not $t$. Strong connectivity therefore implies
$\partial_D^{+}(X_j^{*})=\{t\}$.
If $X_j^{*}\cap O=\varnothing$, every vertex of $X_j^{*}$ has
out-degree two in $D$. Since $g(D)\ge 3k\ge k+1$,
Lemma~\ref{lem:oneexit} would then yield a subdivision of $C(k,k)$,
a contradiction. Thus $X_j^{*}\cap O\ne\varnothing$.
\end{claimproof}

The remaining argument uses the properties of the refined sets and the bounds in~\eqref{eq:finegap}. In particular, it no longer depends on how the original segments were chosen.

\begin{claim}\label{cl:adjacent}
If $X_i^*\cap X_j^*\ne\varnothing$ and $i\ne j$, then $j=i-1$ or $j=i+1$ modulo $N$.
\end{claim}

\begin{claimproof}{cl:adjacent}
Since $X_i^{*}\cap X_j^{*}$ is nonempty and proper, strong connectivity ensures that there is an arc $uv$ with $u\in X_i^{*}\cap X_j^{*}$ and $v\notin X_i^{*}\cap X_j^{*}$.

If $v\notin X_i^{*}\cup X_j^{*}$, the unique external out-neighbor of each set forces $v=\tilde c_{i+1}^{*}=\tilde c_{j+1}^{*}$, contrary to $i\ne j$.

If $v\in X_i^{*}\setminus X_j^{*}$, then $v=\tilde c_{j+1}^{*}$. By~\eqref{eq:fineintersection}, the only selected vertex in $X_i^{*}$ is $\tilde c_i^{*}$. Thus $j+1=i$.

Similarly, if $v\in X_j^{*}\setminus X_i^{*}$, then $i+1=j$. Here all indices are taken modulo $N$.
\end{claimproof}

\begin{claim}\label{cl:hanging}
Suppose that $X_i^{*}\cap X_{i+1}^{*}\ne\varnothing$, and let
$c=\tilde c_{i+1}^{*}$ and $d=\tilde c_{i+2}^{*}$.
Then there is a directed cycle $C'$ meeting $C$ only at $c$.
Moreover, every path from a vertex of $C'-c$ to a vertex of $C-c$
passes through $c$.
\end{claim}

\begin{claimproof}{cl:hanging}
By strong connectivity, there is an arc $uv$ leaving $X_i^*\cap X_{i+1}^*$. If $v$ belongs to neither set, the unique-exit property forces both $v=c$ and $v=d$, a contradiction. If $v\in X_i^*\setminus X_{i+1}^*$, then $v=d$. However, by~\eqref{eq:fineintersection}, the only selected vertex in $X_i^*\cap V(C)=V(C[\tilde c_i^*,c))$ is $\tilde c_i^*$. Since $N\ge5$, the selected vertex $d=\tilde c_{i+2}^*$ is distinct from $\tilde c_i^*$, so $d\notin X_i^*$. Thus $v\in X_{i+1}^*\setminus X_i^*$ and $v=c$. Write the arc as $uc$.

Since $u\in X_{i+1}^*$, choose a $(c,u)$-path $P$ in $D-d$. Every vertex of $P$ lies in $X_{i+1}^*$. Let $x$ be its last vertex on $C$. Then $x\in C[c,d)$. If $x\ne c$, the path $P[x,u]\circ uc$ is internally disjoint from $C$, while
\[
\len{C[x,c]}\ge\len{C[d,c]}
=L-\len{C[c,d]}\ge L-(k-1)\ge k.
\]
This contradicts~\eqref{eq:ear}. Thus $x=c$, and $C'=P\circ uc$ is a directed cycle meeting $C$ only at $c$.

Let $w$ be the first vertex of $C'$ encountered after $c$ that belongs to $X_i^*$. Such a vertex exists because $u\in X_i^*$, and $w\ne c$ because $c\notin X_i^*$. Choose a path $Q$ from $\tilde c_i^*$ to $w$ in $D-c$. Its vertices all lie in $X_i^*$. It follows from~\eqref{eq:fineintersection} that $Q$ meets $C[c,\tilde c_i^*]$ only at $\tilde c_i^*$, and from the choice of $w$ that $Q$ meets $C'[c,w]$ only at $w$. Therefore $C[c,\tilde c_i^*]\circ Q$ and $C'[c,w]$ are internally vertex-disjoint $(c,w)$-paths. The first has length at least $L-(k-1)\ge k$, so the second must have length at most $k-1$. Consequently,
\begin{equation}\label{eq:shortlong}
\len{C'[c,w]}\le k-1\;\text{and}\;
\len{C'[w,c]}\ge g(D)-(k-1)\ge k.
\end{equation}

Figure~\ref{fig:hanging} shows the paths used above. Further intersections of $Q$ with $C'[w,c]$ are not displayed.

\begin{figure}[tbp]
\centering
\begin{tikzpicture}[font=\small]
\node[dot,label=left:{$\tilde c_i^*$}] (a) at (-3,0) {};
\node[dot,label=below:{$c=\tilde c_{i+1}^*$}] (c) at (0,0) {};
\node[dot,label=right:{$d=\tilde c_{i+2}^*$}] (d) at (3,0) {};
\draw[patharc] (a)--node[above=4pt] {$\le k-1$}(c);
\draw[patharc] (c)--node[above=4pt] {$\le k-1$}(d);
\draw[patharc] (d) .. controls (3,-2.5) and (-3,-2.5) .. node[below=3pt] {$C[d,\tilde c_i^*]$} (a);
\node at (0,-1.25) {$C$};
\node[dot,label=above:{$w$}] (w) at (0,3.2) {};
\node[dot,label=right:{$u$}] (u) at (1.25,1.35) {};
\draw[patharc] (c) .. controls (-1.8,1.2) and (-1.3,2.9) .. (w);
\node at (-.15,1.8) {$\le k-1$};
\draw[patharc] (w) .. controls (1.1,3.1) and (1.6,2.3) .. (u);
\draw[realarc] (u)--(c);
\node at (.3,2.5) {$C'$};
\draw[patharc] (a) .. controls (-4.2,2.7) and (-2.1,4.2) .. node[above left=2pt] {$Q$} (w);
\node[align=center] at (4.1,2.0) {$\lVert C'[w,c]\rVert$\\$\ge g(D)-(k-1)$};
\end{tikzpicture}
\caption{The cycles $C,C'$ and the path $Q$.}
\label{fig:hanging}
\end{figure}

Suppose now that a path in $D-c$ joins $C'-c$ to $C-c$. By taking its first vertex on $C-c$ and its last vertex on $C'$ before that vertex, we obtain a path $R$ with initial vertex $y\in C'-c$ and terminal vertex $z\in C-c$, whose internal vertices avoid $C\cup C'$. If $z\ne d$, then $C'[c,y]\circ R$ is a path in $D-d$, so $z\in X_{i+1}^*\cap V(C)$. Equation~\eqref{eq:fineintersection} then implies $z\in C(c,d)$. Thus $z\in C(c,d]$ in either case.

If $y\in C'[w,c)$, the walk $QC'[w,y]R$ avoids $c$ and starts at $\tilde c_i^*$, so $z\in X_i^*$. This again contradicts~\eqref{eq:fineintersection}. Thus $y\in C'(c,w)$. The paths $C'[y,c]$ and $R\circ C[z,c]$ are internally vertex-disjoint. The former contains $C'[w,c]$ and hence has length at least $k$ by~\eqref{eq:shortlong}. The latter contains $C[d,c]$, of length at least $L-(k-1)\ge k$. Their union is the forbidden subdivision, proving the claim.
\end{claimproof}

\begin{claim}\label{cl:concentrate}
If $X_i^*\cap X_{i+1}^*\ne\varnothing$, then $O\subseteq X_i^*\cup X_{i+1}^*$.
\end{claim}

\begin{claimproof}{cl:concentrate}
Use the notation $c,d,C'$ of Claim~\ref{cl:hanging}, and let $q$ and $p$ be the successors of $c$ on $C$ and $C'$, respectively. Since the cycles meet only at $c$, we have $p\ne q$. By~\eqref{eq:minimal}, these are exactly the two out-neighbors of $c$, and $c\notin O$.

Delete the arc $cq$, and define
\[
U=\{v\in V(D):\text{there is a $(\tilde c_i^*,v)$-path in }D-cq\}.
\]
The segment $C[\tilde c_i^*,c]$ avoids $cq$, so $c\in U$. To prove that $U\subseteq X_i^*\cup X_{i+1}^*$, fix $v\in U$ and a simple $(\tilde c_i^*,v)$-path $T$ in $D-cq$.

\proofcase{1}{$c\notin V(T)$.}
Then $T$ is a path in $D-c$, so $v\in X_i^*$ by definition.

\proofcase{2}{$c\in V(T)$.}
If $v=c$, then $v\in X_{i+1}^*$. Otherwise, $T[c,v]$ starts with $cp$, since $cq$ has been deleted and these were the only two arcs leaving $c$. The path $T[p,v]$ avoids $c$ by simplicity. Since $p\in C'-c$, Claim~\ref{cl:hanging} implies that it cannot meet $C-c$: a prefix ending at such a vertex would be a path in $D-c$ from $C'-c$ to $C-c$. In particular it avoids $d$. Thus $T[c,v]$ is a $(c,v)$-path in $D-d$, and $v\in X_{i+1}^*$. We conclude that
\begin{equation}\label{eq:U}
U\subseteq X_i^*\cup X_{i+1}^*.
\end{equation}

Suppose that some $e\in O$ lies outside this union. Then $e\notin U$, so $U$ is a proper subset of $V(D)$. Set $H=(D-cq)[U]$. By the definition of $U$, every out-neighbor in $D-cq$ of a vertex of $U$ also lies in $U$. Hence every vertex of $H$ other than $c$ retains its out-degree from $D$, while the out-degree of $c$ decreases from two to one. It follows that $\delta^+(H)\ge1$,\;$ O(H)=(O\cap U)\cup\{c\}$ and $|O(H)|\le(r-1)+1=r$.
As a subdigraph of $D$, the digraph $H$ has digirth at least $G$ and contains no subdivision of $C(k,k)$. It is nonempty and has fewer vertices than $D$, contradicting the choice of $D$. The claim follows.
\end{claimproof}

\begin{claim}\label{cl:disjoint}
The sets $X_0^*,\ldots,X_{N-1}^*$ are pairwise disjoint.
\end{claim}

\begin{claimproof}{cl:disjoint}
Suppose that two of the sets intersect. By Claim~\ref{cl:adjacent}, they have the form $X_i^*,X_{i+1}^*$. Since $N\ge5$, the set $X_{i+3}^*$ is cyclically nonconsecutive with each of them and is therefore disjoint from their union. Claim~\ref{cl:concentrate} implies $O\subseteq X_i^*\cup X_{i+1}^*$, so $X_{i+3}^*\cap O=\varnothing$, contrary to Claim~\ref{cl:fine}.
\end{claimproof}

By Claim~\ref{cl:fine}, each of the pairwise disjoint sets in Claim~\ref{cl:disjoint} contains a vertex of $O$. Therefore $N\le|O|=r$. Summing~\eqref{eq:finegap} around the cycle yields
\[
g(D)=L=\sum_{j=0}^{N-1}\len{C[\tilde c_j^*,\tilde c_{j+1}^*]}
\le N(k-1)\le r(k-1),
\]
contrary to $g(D)\ge r(k-1)+1$. This proves the upper bound for $k\ge3$.

\subsection{\texorpdfstring{The upper bound for $k=2$}{The upper bound for k = 2}}\label{sec:ktwo}
We now consider $k=2$ and $r\ge4$, with $G=r+1$. Choose a counterexample with $g(D)\ge r+1$ of minimum order and, subject to this, minimum size. The reduction at the start of Section~\ref{sec:upper} ensures strong connectivity, out-degrees in $\{1,2\}$, and exactly $r$ vertices of out-degree one. Write $O=O(D)$. Let $C=x_0x_1\cdots x_{L-1}x_0$ be a shortest directed cycle, where $L\ge r+1\ge5$ and indices are read modulo $L$.

The proof of Claim~\ref{cl:chain} applies without change with $k=2$: consecutive cut-vertices on an isometric segment have distance at most $k-1=1$. Thus, for every arc $uv$ of $C$, every internal vertex of $C[v,u]$ is a $(v,u)$-cut-vertex. Define
\[
\tilde c_i^*=x_i\;\text{and}\;
X_i^*=\{v\in V(D):\text{there is an $(x_i,v)$-path in }D-x_{i+1}\}.
\]
If a path from $x_i$ avoiding $x_{i+1}$ reached another vertex $z$ of $C$, then $z\ne x_{i+1}$, and appending $C[z,x_{i-1}]$ would give an $(x_i,x_{i-1})$-walk avoiding $x_{i+1}$. This contradicts the cut-vertex property for the arc $x_{i-1}x_i$. Consequently $X_i^*\cap V(C)=\{x_i\}$ and $\partial_D^+(X_i^*)=\{x_{i+1}\}$.
By Lemma~\ref{lem:oneexit}, $X_i^*\cap O\ne\varnothing$.

The observation in~\eqref{eq:ear} also holds for $k=2$. Indeed, if $x,y$ are distinct vertices of $C$ and $P$ is an $(x,y)$-path whose internal vertices avoid $C$, then $\len{C[x,y]}\ge2$ would imply $\len P\ge\dist_D(x,y)=\len{C[x,y]}\ge2$ by isometry. Hence $P\cup C[x,y]$ would be a subdivision of $C(2,2)$, a contradiction. Thus $\len{C[x,y]}\le1$. Moreover, $L=g(D)\ge5$ implies $L-1\ge2$ and $g(D)-1\ge2$, so the length estimates used in Claim~\ref{cl:hanging} remain valid.

Here $N=L\ge5$, and consecutive selected vertices are adjacent on $C$. Together with these observations, the proofs of Claims~\ref{cl:adjacent}--\ref{cl:disjoint} apply to these sets, so they are pairwise disjoint. Since each contains a vertex of out-degree one, we obtain $L=N\le r$, contrary to $L\ge r+1$. This completes the proof of the sufficient condition for $k=2$.

\subsection{Sharpness of the bound}\label{sec:lower}
\begin{proposition}\label{prop:construction}
For all integers $k\ge2$ and $r\ge4$, there exists a digraph with exactly $r$ vertices of out-degree one, all remaining vertices of out-degree two, and digirth $r(k-1)$, containing no subdivision of $C(k,k)$.
\end{proposition}

\begin{proof}
Suppose first that $k\ge3$. Take distinct vertices $v_0,\ldots,v_{r-1}$, with indices read modulo $r$. For each $i$, introduce $k-2$ disjoint layers $\{a_{i,j},b_{i,j}\}$ ($1\le j\le k-2$) between $v_i$ and $v_{i+1}$. Add arcs from $v_i$ to both vertices of the first layer, from each vertex of a layer to both vertices of the next layer, and from both vertices of the last layer to $v_{i+1}$. Finally, add the arc $b_{i,k-2}a_{i,k-2}$. There are no other arcs. Figure~\ref{fig:module} illustrates one of these modules when $k=5$.

\begin{figure}[!htb]
\centering
\begin{tikzpicture}[font=\small,scale=.9,transform shape]
\node at (5,-2.6) {(a) Cyclic arrangement of the modules};
\node[gate] (v0) at (0,1) {$v_0$};
\node[gate] (v1) at (3,1) {$v_1$};
\node[gate] (v2) at (6,1) {$v_2$};
\node[gate] (vr) at (10,1) {$v_{r-1}$};
\draw[realarc] (v0)--node[above=4pt] {$M_0$}(v1);
\draw[realarc] (v1)--node[above=4pt] {$M_1$}(v2);
\node at (8,1) {$\cdots$};
\draw[realarc] (v2)--(7.2,1);
\draw[realarc] (8.7,1)--(vr);
\draw[realarc] (vr) .. controls (10,-1.2) and (0,-1.2) .. node[below=3pt] {$M_{r-1}$}(v0);
\node[align=center] at (5,-1.9) {Each arrow in (a) represents a module $M_i$};
\begin{scope}[yshift=-5.6cm]
\node at (5,-3.6) {(b) The module $M_i$ for $k=5$};
\node[gate] (v) at (0,0) {$v_i$};
\node[gate] (w) at (10,0) {$v_{i+1}$};
\foreach \j/\xx in {1/2.3,2/5,3/7.7}{
 \node[dot,label=above:{$a_{i,\j}$}] (a\j) at (\xx,.7) {};
 \node[dot,label=below:{$b_{i,\j}$}] (b\j) at (\xx,-.7) {};
}
\draw[realarc] (v)--(a1);\draw[realarc] (v)--(b1);
\foreach \j/\next in {1/2,2/3}{
 \draw[realarc] (a\j)--(a\next);\draw[realarc] (a\j)--(b\next);
 \draw[realarc] (b\j)--(a\next);\draw[realarc] (b\j)--(b\next);
}
\draw[realarc] (a3)--(w);\draw[realarc] (b3)--(w);
\draw[realarc,line width=1.1pt] (b3)--node[right=4pt] {$e_i$}(a3);
\draw (a3) circle[radius=3.3pt];
\node at (2.3,-1.65) {Layer 1};
\node at (5,-1.65) {Layer 2};
\node at (7.7,-1.65) {Layer 3};
\node[align=center] at (5,-2.65) {$e_i=b_{i,3}a_{i,3}$; the circled $a_{i,3}$ has out-degree one\\Without $e_i$: length $4=k-1$; with $e_i$: length $5=k$};
\end{scope}
\end{tikzpicture}
\caption{The lower-bound construction.}
\label{fig:module}
\end{figure}

The vertices $a_{i,k-2}$ are precisely the vertices of out-degree one; every other vertex has out-degree two. Each module is acyclic, since arcs move to the next layer, except for the single last-layer arc from $b_{i,k-2}$ to $a_{i,k-2}$. Every directed cycle must therefore pass through the vertices $v_i$, in their cyclic order. A passage from $v_i$ to $v_{i+1}$ has length $k-1$ or $k$, according to whether the last-layer arc is used. Hence every cycle has length at least $r(k-1)$, and choosing a passage of length $k-1$ in each module produces a cycle of exactly that length.

Suppose that there are two internally vertex-disjoint $(s,t)$-paths, both of length at least $k$. If $s$ lies inside a module, every path from $s$ that stops before the next vertex $v_{i+1}$ has length at most $k-2$, and every path that first reaches $v_{i+1}$ has length at most $k-1$. Each of the two long paths must therefore contain $v_{i+1}$ internally, a contradiction.

We may consequently assume that $s=v_i$. A path starting at $v_i$ and avoiding $v_{i+1}$ has length at most $k-1$, so both paths pass through $v_{i+1}$. If $t\ne v_{i+1}$, they have this vertex in common internally. If $t=v_{i+1}$, simplicity forces both paths to remain within the module until their terminal vertex. For either path to have length at least $k$, it must use the arc $b_{i,k-2}a_{i,k-2}$. Thus they again share internal vertices. This excludes a subdivision of $C(k,k)$.

For $k=2$, take a directed cycle of length $r$. Every vertex has out-degree one, its digirth is $r$, and between any two distinct vertices there is only one directed path. Hence it contains no subdivision of $C(2,2)$.
\end{proof}

Proposition~\ref{prop:construction} implies $F_r(k)\ge r(k-1)+1$, while Sections~\ref{sec:upper} and~\ref{sec:ktwo} prove the opposite inequality. This completes the proof of Theorem~\ref{thm:many}.

\section{Concluding remarks}\label{sec:remarks}
Theorem~\ref{thm:many} determines $F_r(k)$ for every $k\ge2$ and $r\ge4$. Together with the result $F_0(k)=F_1(k)=k+1$ of Chen et al.~\cite{chenexact}, this leaves the cases of two and three permitted exceptions.

\begin{problem}\label{prob:remaining}
Determine $F_2(k)$ and $F_3(k)$ for every integer $k\ge2$.
\end{problem}

The distinction between one exception and several exceptions is reflected in the proof. The local lemma forces every region with a unique external out-neighbor to contain an exceptional vertex. The argument for several exceptions requires, in addition, control of intersections between such regions. The resulting disjointness permits a direct count along a shortest directed cycle and yields the sharp bound $r(k-1)+1$.

\section*{Declaration on the use of AI}
The authors used AI tools in preparing this manuscript to explore proof strategies, check mathematical arguments, improve the exposition, and prepare LaTeX text and TikZ illustrations. The authors have reviewed and verified all proofs, arguments,
and figures, and take full responsibility for the content of this publication.

\end{document}